\documentstyle{article}
\begin{document}

\def\dinf{\displaystyle\inf}
\def\dsup{\displaystyle\sup}
\def \theequation {\arabic {equation}}
\setcounter {equation}{0}

\vspace*{1cm}
\begin{center}
{{{\LARGE {\bf Robust Strictly Positive Real Synthesis for
Polynomial Families of Arbitrary Order} \footnote{Supported by
National Natural Science Foundation of China (No. 60204006 and No.
69925307). }}}}
\end{center}

\vskip 0.6cm \centerline{ WenshengYu \ LongWang } \vskip 6pt

\centerline{\small {\it Center for Systems and Control }}

\centerline{\small {\it Department of Mechanics and Engineering
Science, Peking University}}

\centerline{\small\it {Beijing  100871, P. R. China, E-mail:
longwang@mech.pku.edu.cn}}

\vskip 0.6cm \noindent { {{Abstract\quad}}} {{{For any two
$n$-$th$ order polynomials $a(s)$ and $b(s),$ the Hurwitz
stability of their convex combination is necessary and sufficient
for the existence of a polynomial $c(s)$ such that $c(s)/a(s)$ and
$c(s)/b(s)$ are both strictly positive real. }}}

\noindent {{ {Keywords\quad}}} {{{Robust Stability, Strict
Positive Realness, Robust Analysis and Synthesis, Polynomial
Segment}}

\strut

\section{Introduction}

{ { The strict positive realness (SPR) of transfer functions is an
important performance specification, and plays a critical role in
various fields such as absolute stability/hyperstability theory
\cite{Kal63,Pop73}, passivity analysis \cite{DV75}, quadratic
optimal control \cite{AM70} and adaptive system theory
\cite{Lan79}. In recent years, stimulated by the parametrization
method in robust stability analysis \cite{ABK93,Bar94,BCK95}, the
study of robust strictly positive real systems has received much
attention, and great progress has beem made \cite{ADK90,BZ93}
\cite{Bia02}-\cite{DB87} \cite{Hen02}-\cite{HHX90}
\cite{MA98}-\cite{PD97} \cite{WH91}-\cite{YWT99}. However, most
results belong to the category of robust SPR analysis. Valuable
results in robust SPR synthesis are rare. The following
fundamental problem is still open \cite{HHX90,MP01,WY99,WY00}
\cite{WY01a}-\cite{YH99} \cite{YW00} \cite{YW01}-\cite{YW03}:

\textit{Suppose $a(s)$ and $b(s)$ are two $n$-$th$ order Hurwitz
polynomials, does there exist, and how to find a (fixed)
polynomial $c(s)$ such that $c(s)/a(s)$ and $c(s)/b(s)$ are both
SPR?}

By the definition of SPR, it is easy to know that the Hurwitz
stability of the convex combination of $a(s)$ and $b(s)$ is
necessary for the existence of polynomial $c(s)$ such that
$c(s)/a(s)$ and $c(s)/b(s)$ are both SPR. In
\cite{HHX89,HHX90,PD97}, it was proved that, if $a(s)$ and $b(s)$
have the same even (or odd) parts, such a polynomial $c(s)$ always
exists; In \cite{ADK90, HHX89,HHX90,MA98,
WY99,WY00,WY01a,Yu98,YW00,YW01}, it was proved that, if $n\leq 4$
and $a(s),b(s)\in K$ ($K$ is a stable interval polynomial set),
such a polynomial $c(s)$ always exists; Recent results show that
\cite{WY99,WY00,WY01a,Yu98,YH99,YW00} \cite{YW01a}-\cite{YW03}, if
$n\leq 6$ and $a(s)$ and $b(s)$ are the two endpoints of the
convex combination of stable polynomials, such a polynomial $c(s)$
always exists. Some sufficient condition for robust SPR synthesis
are presented in \cite{ADK90,BZ93,DB87,MA98,WY99,WY00,Yu98},
especially, the design method in \cite{WY99,WY00} is numerically
efficient for high-order polynomial segments or interval
polynomials, and the derived conditions are necessary and
sufficient for low-order polynomial segments or interval
polynomials.

This paper presents a constructive proof of the following
statement: for any two $n$-$th$ order polynomials $a(s)$ and
$b(s),$ the Hurwitz stability of their convex combination is
necessary and sufficient for the existence of a polynomial $c(s)$
such that $c(s)/a(s)$ and $c(s)/b(s)$ are both SPR. This also
shows that the conditions given in \cite{WY99,WY00} are also
necessary and sufficient and the above open problem admits a
positive answer. Some previously obtained SPR synthesis results
for low-order polynomial segments
\cite{WY99,WY00,WY01a,Yu98,YH99,YW00} \cite{YW01a}-\cite{YW03}
become special cases of our main result in this paper.

\section{Main Results}

In this paper, $ P^n$ stands for the set of $n$-$th$ order
polynomials of $s$ with real coefficients, $ R$ stands for the
field of real numbers, $ \partial (p)$ stands for the order of
polynomial $ p(\cdot )$, and $H^n\subset P^n$ stands for the set
of $n$-$th$ order Hurwitz stable polynomials with real
coefficients.

In the following definition, $ p(\cdot )\in P^m, q(\cdot )\in P^n, f(s)=p(s)/q(s)$ is
a rational function.

{\bf  Definition 1} \cite{YW99} \ \ $ f(s)$ is said to be strictly
positive real(SPR), if

(i) $ \partial (p)=\partial (q);$

(ii) $f(s)$ is analytic in $ \mbox{Re}  [s]\geq 0,$  (namely, $ q(\cdot )\in H^n$ );

(iii) $ \mbox{Re}  [f(j\omega )]>0,\ \ \ \forall \omega \in R.$

If $f(s)=p(s)/q(s)$ is proper, it is easy to get the following property:

{\bf Property 1} \cite{CDB91} \ \ If $f(s)=p(s)/q(s)$ is a proper
rational function, $q(s)\in H^n,$ and $\forall \omega \in R,
\mbox{Re}  [f(j\omega )]>0,$ then $p(s)\in H^n\cup H^{n-1}.$

The following theorem is the main result of this paper:

{\bf  Theorem 1} \ \ Suppose $ a(s)=s^n+a_1s^{n-1}+\cdots+a_n\in
H^n, b(s)=s^n+b_1s^{n-1}+\cdots+b_n\in H^n,$ the necessary and
sufficient condition for the existence of a polynomial $ c(s)$
such that $ c(s)/a(s)$ and $ c(s)/b(s)$ are both Strict Positive
Real, is
$$ \lambda b(s)+(1-\lambda )a(s)\in H^n,\lambda \in [0,1].$$

The statement is obviously true for the cases when $n=1$ or $n=2.$
We will prove it for the case when $n\geq3.$

Since SPR transfer functions enjoy convexity property, by Property
1, we get the necessary part of the theorem.

To prove sufficiency, we must first introduce some lemmas.

{\bf Lemma {\bf 1}}\ \ Suppose $ a(s)=s^n+a_1s^{n-1}+\cdots+a_n\in
H^n,$ then, for every $k\in\{1,2,\cdots,n-2\},$ the following
quadratic curve is an ellipse in the first quadrant (i.e.,
$x_i\geq 0, i=1,2,\cdots, n-1$) of the $R^{n-1}$ space
$(x_1,x_2,\cdots,x_{n-1})$
\footnote{ When $n=3,$ the ellipse
equation is (see \cite{Yu98, YH99} for details):
$$
(a_2x_1-a_1x_2-a_3)^2-4(a_1-x_1)(a_3x_2)=0.
$$

When $n=4,$ the two ellipse equations are (see \cite{WY01a}
\cite{YW00}-\cite{YW01a} for details): $$
\begin{array}{l}
\left\{
\begin{array}{l}
(a_2x_1+x_3-a_1x_2-a_3)^2-4(a_1-x_1)(a_3x_2-a_2x_3-a_4x_1)=0, \\
a_4x_3=0,
\end{array}
\right.  \\
\left\{
\begin{array}{l}
(a_3x_2-a_2x_3-a_4x_1)^2-4(a_2x_1+x_3-a_1x_2-a_3)a_4x_3=0, \\
a_1-x_1=0.
\end{array}
\right.  \\
\end{array}
$$

When $n=5,$ the three ellipse equations are (see \cite{YW01b} for
details): $$
\begin{array}{l}
\left\{
\begin{array}{l}
(a_2x_1+x_3-a_1x_2-a_3)^2-4(a_1-x_1)(a_5+a_3x_2+a_1x_4-a_2x_3-a_4x_1)=0, \\
a_4x_3-a_3x_4-a_5x_2=0, \ a_5x_4=0,\\
\end{array}
\right. \\
\left\{
\begin{array}{l}
(a_5+a_3x_2+a_1x_4-a_2x_3-a_4x_1)^2-4(a_2x_1+x_3-a_1x_2-a_3)(a_4x_3-a_3x_4-a_5x_2)=0,
\\
a_1-x_1=0, \ a_5x_4=0,\\
\end{array}
\right. \\
\left\{
\begin{array}{l}
(a_4x_3-a_3x_4-a_5x_2)^2-4(a_5+a_3x_2+a_1x_4-a_2x_3-a_4x_1)a_5x_4=0, \\
a_1-x_1=0, \ a_2x_1+x_3-a_1x_2-a_3=0.\\
\end{array}
\right. \\
\end{array}
$$

When $n=6,$ the four ellipse equations are (see \cite{YW03} for
details): $$
\begin{array}{l}
\left\{
\begin{array}{l}
(a_2x_1+x_3-a_1x_2-a_3)^2-4(a_1-x_1)(a_5+a_3x_2+a_1x_4-x_5-a_2x_3-a_4x_1)=0,
\\
a_6x_1+a_4x_3+a_2x_5-a_3x_4-a_5x_2=0, \ a_5x_4-a_4x_5-a_6x_3=0, \ a_6x_5=0,\\
\end{array}
\right. \\
\left\{
\begin{array}{l}
(a_5+a_3x_2+a_1x_4-x_5-a_2x_3-a_4x_1)^2-4(a_2x_1+x_3-a_1x_2-a_3)(a_6x_1+a_4x_3+a_2x_5-a_3x_4-a_5x_2)=0,
\\
a_1-x_1=0, \ a_5x_4-a_4x_5-a_6x_3=0, \ a_6x_5=0,\\
\end{array}
\right. \\
\left\{
\begin{array}{l}
(a_6x_1+a_4x_3+a_2x_5-a_3x_4-a_5x_2)^2-4(a_5+a_3x_2+a_1x_4-x_5-a_2x_3-a_4x_1)(a_5x_4-a_4x_5-a_6x_3)=0
\\
a_1-x_1=0, \ a_2x_1+x_3-a_1x_2-a_3=0, \ a_6x_5=0, \\
\end{array}
\right. \\
\left\{
\begin{array}{l}
(a_5x_4-a_4x_5-a_6x_3)^2-4(a_6x_1+a_4x_3+a_2x_5-a_3x_4-a_5x_2)a_6x_5=0, \\
a_1-x_1=0, \ a_2x_1+x_3-a_1x_2-a_3=0, \ a_5+a_3x_2+a_1x_4-x_5-a_2x_3-a_4x_1=0.\\
\end{array}
\right. \\
\end{array}
$$
}:
\begin{center}
$\left\{
\begin{array}{l}
c_{k+1}^2-4c_{k}c_{k+2}=0, \\
c_l=0, \\
l\in \{1,2,\cdots,n\},l\neq k,k+1,k+2,
\end{array}
\right. $
\end{center}
and this ellipse is tangent with the line
\begin{center}
$\left\{
\begin{array}{l}
c_l=0, \\
l\in \{1,2,\cdots,n\},l\neq k+1,k+2,
\end{array}
\right. $
\end{center}
and the line
\begin{center}
$\left\{
\begin{array}{l}
c_l=0,  \\
l\in \{1,2,\cdots,n\},l\neq k,k+1,
\end{array} \right. $
\end{center}
respectively, where
$c_l:=\sum\limits_{j=0}^{n}(-1)^{l+j}a_{j}x_{2l-j-1},l=1,2,\cdots,n,$
 $a_0=1,x_0=1,$ $a_i=0$ if $i<0$ or $i>n,$ and $x_i=0$ if
$i<0$ or $i>n-1.$

{\bf Proof}\ \ Since $ a(s)$ is Hurwitz stable, by using
mathematical induction, Lemma 1 is proved by a direct calculation.

For notational simplicity, for $ a(s)=s^n+a_1s^{n-1}+\cdots+a_n\in
H^n, b(s)=s^n+b_1s^{n-1}+\cdots+b_n\in H^n,$ $\forall
k\in\{1,2,\cdots,n-2\},$ denote
$$\begin{array}{ll}
\Omega _{ek}^a:=\{(x_1,x_2,\cdots,x_{n-1})|&
c_{k+1}^2-4c_{k}c_{k+2}<0, \\
&c_l=0,l\in \{1,2,\cdots,n\},l\neq k,k+1,k+2 \},
\end{array}
$$ %
and
$$\begin{array}{ll}
\Omega _{ek}^b:=\{(x_1,x_2,\cdots,x_{n-1})|&
d_{k+1}^2-4d_{k}d_{k+2}<0, \\
&d_l=0,l\in \{1,2,\cdots,n\},l\neq k,k+1,k+2 \},
\end{array}
$$ %
where $c_l:=\sum\limits_{j=0}^{n}(-1)^{l+j}a_{j}x_{2l-j-1},
d_l:=\sum\limits_{j=0}^{n}(-1)^{l+j}b_{j}x_{2l-j-1},l=1,2,\cdots,n,$
 $a_0=1,b_0=1,x_0=1,$ $a_i=0$ and $b_i=0$ if $i<0$ or $i>n,$
and $x_i=0$ if $i<0$ or $i>n-1.$

In what follows, $(A, B)$ stands for the set of points in the line
segment connecting the point $A$ and the point $B$ in the
$R^{n-1}$ space $(x_1,x_2,\cdots,x_{n-1}),$ not including the
endpoints $A$ and $B$. Denote
$$\begin{array}{ll}
\Omega ^a:=\{(x_1,x_2,\cdots,x_{n-1})|
&(x_1,x_2,\cdots,x_{n-1})\in \bigcup_{i=1,i<j\leq
n-2}^{n-3}(A_i,A_j),
\\
& \forall A_i \in \Omega _{ei}^a,i\in \{1,2,\cdots,n-2\} \}
\end{array}$$
and
$$\begin{array}{ll}
\Omega ^b:=\{(x_1,x_2,\cdots,x_{n-1})|
&(x_1,x_2,\cdots,x_{n-1})\in \bigcup_{i=1,i<j\leq
n-2}^{n-3}(B_i,B_j),
\\
& \forall B_i \in \Omega _{ei}^b,i\in \{1,2,\cdots,n-2\} \}.
\end{array}$$

{\bf  Lemma {\bf 2}}\ \ Suppose $
a(s)=s^n+a_1s^{n-1}+\cdots+a_n\in H^n,
b(s)=s^n+b_1s^{n-1}+\cdots+b_n\in H^n,$ if $ \Omega ^a \cap \Omega
^b\neq \phi,$ take $(x_1,x_2,\cdots,x_{n-1})\in \Omega ^a \cap
\Omega ^b, $ and let $
c(s):=s^{n-1}+(x_1-\varepsilon)s^{n-2}+x_2s^{n-3}+\cdots+x_{n-2}s+(x_{n-1}+\varepsilon)$
($ \varepsilon $ is a sufficiently small positive number), then
for $\displaystyle\frac{c(s)}{a(s)}$ and
$\displaystyle\frac{c(s)}{b(s)}$, we have $\forall \omega \in R,
\mbox{Re} [\displaystyle\frac {c(j\omega )}{a(j\omega )}]>0$  and
$ \mbox{Re} [\displaystyle\frac {c(j\omega )}{b(j\omega )}]>0.$

{\bf Proof}\ \ Suppose $(x_1,x_2,\cdots,x_{n-1})\in \Omega ^a \cap \Omega ^b, $ let $%
c(s):=s^{n-1}+(x_1-\varepsilon)s^{n-2}+x_2s^{n-3}+\cdots+x_{n-2}s+(x_{n-1}+\varepsilon),
\varepsilon
>0,\varepsilon $ sufficiently small.

$\forall \omega \in R,$ consider
$$
\begin{array}{lll}
\mbox{Re}[\displaystyle\frac{c(j\omega )}{a(j\omega )}] & = & \displaystyle%
\frac 1{\mid a(j\omega )\mid ^2}[c_1\omega ^{2(n-1)}+c_2\omega
^{2(n-2)}+\cdots+c_{n-1}\omega ^2+c_n]\\
& &+\mbox{Re}[\displaystyle\frac{-\varepsilon(j\omega
)^{n-2}+\varepsilon}{a(j\omega )}] \\
& = & \displaystyle%
\frac 1{\mid a(j\omega )\mid ^2}[c_1\omega ^{2(n-1)}+c_2\omega
^{2(n-2)}+\cdots+c_{n-1}\omega ^2+c_n]\\
& &+\displaystyle \frac{(-\varepsilon )}{|a(j\omega )|^2}(-\omega
^{2(n-1)}+\tilde c(\omega^2 )),
\end{array}
$$
where
$c_l:=\sum\limits_{j=0}^{n}(-1)^{l+j}a_{j}x_{2l-j-1},l=1,2,\cdots,n,$
$a_0=1,x_0=1,$ $a_i=0$ if $i<0$ or $i>n,$ and $x_i=0$ if $i<0$ or
$i>n-1,$ and $\tilde c(\omega^2 )$ is a real polynomial with order
not greater than $2(n-2)$.

In order to prove that $\forall \omega \in R,\mbox{Re}[\displaystyle\frac{%
c(j\omega )}{a(j\omega )}]>0,$ let $t=\omega ^2,$ we only need to
prove that, for any $\varepsilon >0,\varepsilon $ sufficiently
small, the following polynomial $f_1(t)$ satisfies
$$
\begin{array}{ll}
f_1(t):= & c_1t ^{n-1}+c_2t^{n-2}+\cdots+c_{n-1}t+c_n \\
& +\varepsilon (t^{n-1}-\tilde c(t ))>0, \ \ \forall t\in
[0,+\infty ).
\end{array}
$$

Since $(x_1,x_2,\cdots,x_{n-1})\in \Omega ^a,$ by the definition
of $\Omega ^a,$ it is easy to know that
$$
g_1(t):= c_1t ^{n-1}+c_2t^{n-2}+\cdots+c_{n-1}t+c_n>0,\ \ \forall
t\in (0,+\infty ).
$$
Moreover, we obviously have $f_1 (0)>0,$ and for any $\varepsilon >0$, when $%
t$ is a sufficiently large or sufficiently small positive number, we have $%
f_1 (t)>0,$ namely, there exist $0<t_1<t_2$ such that, for all
$\varepsilon
>0$, $t\in [0,t_1]\cup [t_2,+\infty)$, we have $f_1 (t)>0.$

Denote
$$
M_1=\dinf_{t\in [t_1,t_2]}g_1(t),
$$
$$
N_1=\dsup_{t\in [t_1,t_2]}|t^{n-1}-\tilde c(t ))|.
$$
Then $M_1>0$ and $N_1>0.$ Choosing $0<\varepsilon <\displaystyle\frac{M_1}{%
N_1},$ by a direct calculation, we have
$$
\begin{array}{ll}
f_1(t):= & c_1t ^{n-1}+c_2t^{n-2}+\cdots+c_{n-1}t+c_n \\
& +\varepsilon (t^{n-1}-\tilde c(t ))>0, \ \ \forall t\in
[0,+\infty ).
\end{array}
$$
Namely,
$$
\forall \omega \in R,\mbox{Re}[\displaystyle\frac{c(j\omega )}{a(j\omega )}%
]>0.
$$

Similarly, since $(x_1,x_2,\cdots,x_{n-1})\in \Omega ^b,$ there
exist $0<t_3<t_4$ such
that, for all $\varepsilon >0$, $t\in [0,t_3]\cup [t_4,+\infty )$, we have $%
f_2(t)>0,$where
$$
\begin{array}{ll}
f_2(t):= & d_1t ^{n-1}+d_2t^{n-2}+\cdots+d_{n-1}t+d_n \\
& +\varepsilon (t^{n-1}-\tilde d(t )),
\end{array}
$$
$$
d_l:=\sum\limits_{j=0}^{n}(-1)^{l+j}b_{j}x_{2l-j-1},\ \
l=1,2,\cdots,n,
$$
where $b_0=1,x_0=1,$ $b_i=0$ if $i<0$ or $i>n,$ and $x_i=0$ if
$i<0$ or $i>n-1,$ and $\tilde d(\omega^2 )$ is a real polynomial
with order not greater than $2(n-2)$ which is determined by the
following equation:
$$
\mbox{Re}[\displaystyle\frac{-\varepsilon(j\omega
)^{n-2}+\varepsilon}{b(j\omega )}]=\displaystyle
\frac{(-\varepsilon )}{|b(j\omega )|^2}(-\omega ^{2(n-1)}+\tilde
d(\omega^2 )).
$$

Denote
$$
g_2(t):= d_1t ^{n-1}+d_2t^{n-2}+\cdots+d_{n-1}t+d_n,
$$
$$
M_2=\dinf_{t\in [t_3,t_4]}g_2(t),
$$
$$
N_2=\dsup_{t\in [t_3,t_4]}|t^{n-1}-\tilde d(t ))|.
$$
Then $M_2>0$ and $N_2>0.$ Choosing $0<\varepsilon <\displaystyle\frac{M_2}{%
N_2},$ we have
$$
\forall \omega \in R,\mbox{Re}[\displaystyle\frac{c(j\omega )}{b(j\omega )}%
]>0.
$$
Thus, by choosing $0<\varepsilon <\min \{{\displaystyle\frac{M_1}{N_1},%
\displaystyle\frac{M_2}{N_2}}\},$ Lemma 2 is proved.

{\bf  Lemma {\bf 3}}\ \ Suppose $
a(s)=s^n+a_1s^{n-1}+\cdots+a_n\in H^n,
b(s)=s^n+b_1s^{n-1}+\cdots+b_n\in H^n,$ if $\lambda
b(s)+(1-\lambda )a(s)\in H^n,\lambda \in [0,1],$ then $ \Omega ^a
\cap \Omega ^b\neq \phi $

{\bf  Proof }\ \  If $ \forall \lambda \in
[0,1],{a_\lambda}(s):=\lambda b(s)+(1-\lambda )a(s)\in H^n,$ by
Lemma 1, for any $ \lambda \in [0,1],$ $\Omega _{ek}^{a_\lambda},
k=1,2,\cdots, n-2,$ are $n-2$ ellipses in the first quadrant of
the $R^{n-1}$ space $(x_1,x_2,\cdots,x_{n-1})$.

$ \forall \lambda \in [0,1],$ denote
$$
\begin{array}{ll}
\Omega ^{a_\lambda}:=\{(x_1,x_2,\cdots,x_{n-1})|
&(x_1,x_2,\cdots,x_{n-1})\in \bigcup_{i=1, i<j\leq
n-2}^{n-3}(A_{\lambda i},A_{\lambda j}),
\\
& \forall A_{\lambda i} \in \Omega _{ei}^{a_\lambda},i\in
\{1,2,\cdots,n-2\} \}
\end{array}
$$ %

Apparently, when $\lambda $ changes continuously from $0$ to $1$,
$\Omega^{a_\lambda}$ will change continuously from $\Omega^{a}$ to
$\Omega^{b},$ and  $\Omega_{ek}^{a_\lambda}$ will change
continuously from $\Omega_{ek}^{a}$ to $\Omega_{ek}^{b},
k=1,2,\cdots,n-2.$

Now assume $ \Omega ^a \cap \Omega ^b= \phi,$ by the definitions
of $ \Omega ^a $ and $ \Omega ^b$ and Lemma 1, $ \exists u_1
> 0,$ $u_2 > 0,u_1\neq a_1, u_1\neq b_1, $ and $\exists \tilde k\in \{1,2,\cdots,n-2\},$
such that the following hyperplane $L$ in the $R^{n-1}$ space
$(x_1,x_2,\cdots,x_{n-1})$
$$
L:\ \ \  \displaystyle\frac  {x_1}{u_1}+\displaystyle\frac
{x_2}{u_2} +\cdots+\displaystyle\frac  {x_{n-1}}{u_{n-1}}=1
$$
separates $\Omega^{a}$ and $\Omega^{b},$ meanwhile, $L$ is tangent
with $ \Omega _{e1}^a ,\Omega _{e2}^a, \cdots,\Omega _{e(n-2)}^a$
and $ \Omega _{e \tilde k}^b $ simultaneously (or tangent with $
\Omega _{e1}^b ,\Omega _{e2}^b, \cdots,\Omega _{e(n-2)}^b$ and $
\Omega _{e \tilde k}^a $ simultaneously).%

Without loss of generality, suppose that $L$ is tangent with $
\Omega _{e1}^a ,\Omega _{e2}^a, \cdots,\Omega _{e(n-2)}^a$ and $
\Omega _{e \tilde k}^b $ simultaneously.%

In what follows, the notation $[x]$ stands for the largest integer
that is smaller than or equal to the real number $x$, and
$\left\langle y\right\rangle _z$ stands for the remainder of the
nonnegative integer $y$ divided by the positive integer $z$
\footnote{ For example, $[1.5]=1,[0.5]=0,[-1.5]=-2$, and
$\left\langle 0\right\rangle _2=0,\left\langle 1\right\rangle
_2=1,\left\langle 11\right\rangle _3=2$.}.

Since $ L$ is tangent with $ \Omega _{e1}^a ,\Omega _{e2}^a,
\cdots,\Omega _{e(n-2)}^a$ and $ \Omega _{e \tilde k}^b $
simultaneously, note that $ a(s)$ is Hurwitz stable and $u_1 >
0,u_1\neq a_1, u_2 > 0, $ using mathematical induction, by a
lengthy calculation, we know that the necessary and sufficient
condition for $ L$ being tangent with $ \Omega _{e1}^a ,\Omega
_{e2}^a, \cdots,\Omega _{e(n-2)}^a$ simultaneously is \footnote{
When $n=3,$ we have (see \cite{Yu98, YH99} for details):
$$
u_1u_2-a_1u_2-a_2u_1+a_3=0.
$$

When $n=4,$ we have (see \cite{WY01a} \cite{YW00}-\cite{YW01a} for
details):
$$
\left\{
\begin{array}{l}
u_1u_2^2-a_1u_2^2-a_2u_1u_2+a_3u_2+a_4u_1=0, \\
u_3=-u_1u_2.
\end{array}
\right.
$$

When $n=5,$ we have (see \cite{YW01b} for details):
$$
\left\{
\begin{array}{l}
u_1u_2^2-a_1u_2^2-a_2u_1u_2+a_3u_2+a_4u_1-a_5=0, \\
u_3=-u_1u_2,\ u_4=-u_2^2.
\end{array}
\right.
$$

When $n=6,$ we have (see \cite{YW03} for details):
$$
\left\{
\begin{array}{l}
u_1u_2^3-a_1u_2^3-a_2u_1u_2^2+a_3u_2^2+a_4u_1u_2-a_5u_2-a_6u_1=0, \\
u_3=-u_1u_2,\ u_4=-u_2^2,\ u_5=u_1u_2^2.
\end{array}
\right.
$$ }
\begin{equation}
\sum_{i=0}^n(-1)^{\displaystyle[\frac{i+1}2]}a_iu_1^{\left\langle
i+1\right\rangle _2}u_2^{\displaystyle[\frac
n2]-\displaystyle[\frac i2]}=0 \label{eq1}
\end{equation}
and \begin{equation}
u_j=(-1)^{\displaystyle[\frac{j-1}2]}u_1^{\left\langle
j\right\rangle _2}u_2^{\displaystyle[\frac{j}2]},\ \
j=3,4,\cdots,n-1, \label{eq2}
\end {equation}
where $a_0=1.$

Since $u_j=(-1)^{\displaystyle[\frac{j-1}2]}u_1^{\left\langle
j\right\rangle _2}u_2^{\displaystyle[\frac{j}2]},\ \
j=3,4,\cdots,n-1, L$ is tangent with $ \Omega _{e \tilde k}^b $,
by a direct calculation, we have
\begin{equation}
\sum_{i=0}^n(-1)^{\displaystyle[\frac{i+1}2]}b_iu_1^{\left\langle
i+1\right\rangle _2}u_2^{\displaystyle[\frac
n2]-\displaystyle[\frac i2]}=0 \label{eq3}
\end{equation}
where $b_0=1.$

From $(\ref{eq1})$  and $(\ref{eq3}),$ we obviously have $\forall
\lambda \in [0,1],$
\begin{equation}
\sum_{i=0}^n(-1)^{\displaystyle[\frac{i+1}2]}a_{\lambda
i}u_1^{\left\langle i+1\right\rangle _2}u_2^{\displaystyle[\frac
n2]-\displaystyle[\frac i2]}=0 \label{eq4}
\end{equation}\\
where $a_{\lambda i}:=a_i+\lambda
(b_i-a_i),a_0=1,b_0=1,i=0,1,2,\cdots,n.$ $(\ref{eq4})$  and
$(\ref{eq2})$ show that $ L$ is also tangent with $ \Omega _{e
\tilde k}^{a_{\lambda }}( \forall \lambda \in [0,1])$, but $L$
separates $ \Omega _{e \tilde k}^a $ and $ \Omega _{e \tilde k}^b,
$ and when $\lambda $ changes continuously from $0$ to $1$,
$\Omega_{e \tilde k}^{a_\lambda}$ will change continuously from
$\Omega_{e \tilde k}^{a}$ to $\Omega_{e \tilde k}^{b},$ which is
obviously impossible. This completes the proof.


{\bf  Lemma {\bf 4}}\ \ Suppose $
a(s)=s^n+a_1s^{n-1}+\cdots+a_n\in H^n,
b(s)=s^n+b_1s^{n-1}+\cdots+b_n\in H^n,
c(s)=s^{n-1}+x_1s^{n-2}+\cdots+x_{n-1},$ if $ \forall \omega \in
R, \mbox{Re} [\displaystyle\frac {c(j\omega )}{a(j\omega )}]>0$
and $ \mbox{Re} [\displaystyle\frac {c(j\omega )}{b(j\omega
)}]>0,$ take
\begin{center}
$ \stackrel{\sim }{c}(s):=c(s)+\delta \cdot h(s),\ \ \delta
>0,\ \delta $ sufficiently small,
\end{center}
(where $h(s)$ is an arbitrarily given monic $n$-$th$ order
polynomial), then $ \displaystyle\frac {\stackrel{\sim
}{c}(s)}{a(s)}$  and
$ \displaystyle\frac {\stackrel{\sim }{c}(s)}{%
b(s)}$ are both strictly positive real.

{\bf  Proof }\ \ Obviously,  $ \partial (\tilde {c})=\partial
(a)=n,$ namely, $ \tilde {c}(s)$ and $ a(s)$ have the same order
$n$. Since $a(s)\in H^n,$ there exists $ \omega _1 >0$ such that,
for all $ \mid \omega \mid \geq \omega _1,$ we have Re$
(\displaystyle\frac {\tilde {c}(j\omega )}{a(j\omega )})>0.$

Denote
$$
M_3=\dinf_{\mid \omega \mid \leq \omega _1}\mbox{Re} (
\displaystyle\frac {c(j\omega )}{a(j\omega )}),  \ \  \ \
N_3=\dsup_{\mid \omega \mid \leq \omega _1}\displaystyle {\mid
\mbox{Re}(\displaystyle\frac {h(j\omega )}{a(j\omega )}) \mid }.
$$
Then $M_3>0 $ and $N_3>0.$ Choosing $0<\varepsilon
<\displaystyle\frac {M_3}{N_3},$ it can be directly verified that
$$
\mbox{Re} (\displaystyle\frac {\tilde {c}(j\omega )}{a(j\omega
)})>0, \forall \omega \in R.
$$

Similarly, $ \partial (\tilde {c})=\partial (b)=n,$ and there
exists $ \omega _2
>0$ such that, for all $ \mid \omega \mid \geq \omega _2,$ we have
Re$ (\displaystyle\frac {\tilde {c}(j\omega )}{b(j\omega )})>0.$

Denote
$$
M_4=\dinf_{\mid \omega \mid \leq \omega _2}\mbox{Re} (
\displaystyle\frac {c(j\omega )}{b(j\omega )}),  \ \  \ \
N_4=\dsup_{\mid \omega \mid \leq \omega _2}\displaystyle {\mid
\mbox{Re}(\displaystyle\frac {h(j\omega )}{b(j\omega )}) \mid }.
$$
Then $M_4>0 $ and $N_4>0.$ Choosing $0<\varepsilon
<\displaystyle\frac {M_4}{N_4},$ it can be directly verified that
$$
\mbox{Re} (\displaystyle\frac {\tilde {c}(j\omega )}{b(j\omega
)})>0, \forall \omega \in R.
$$
Thus, by choosing $0<\varepsilon <\min \{{\displaystyle\frac{M_3}{N_3},%
\displaystyle\frac{M_4}{N_4}}\},$ Lemma 4 is proved.

The sufficiency of Theorem 1 is now proved by simply combining Lemmas 1-4. %

{\bf  Remark {\bf 1}} \ \ From the proof of Theorem 1, we can see
that this paper not only proves the existence, but also provides a
design method. In fact, based on the main idea of this paper, we
have developed a geometric algorithm with order reduction for
robust SPR synthesis which is very efficient
for high order polynomial segments \cite{XWY02}.%

{\bf  Remark {\bf 2}} \ \ The method provided in this paper is
constructive, and is insightful and helpful in solving more
general
robust SPR synthesis problems for polynomial polytopes, multilinear families, etc.. %

{\bf  Remark {\bf 3}} \ \ Our main results in this paper can also
be extended to discrete-time case. In fact, by applying the
bilinear transformation, we can transform the unit circle into the
left half plane. Hence, Theorem 1 can be generalized to
discrete-time case. Moreover, in discrete-time case, the order of
the polynomial obtained by our method is bounded by the order
of the given polynomial segment \cite{Yu98,YW99}.%

{\bf  Remark {\bf 4}} \ \
If $\displaystyle\frac{c(s)}{a(s)}$ and $\displaystyle\frac{c(s)}{b(s)}$ are both SPR,
it is easy to know that
$\forall \lambda \in [0, 1], \displaystyle \frac{c(s)}{\lambda a(s) + (1- \lambda) b(s)}$ is also SPR.

{\bf  Remark {\bf 5}} \ \ The stability of polynomial segment can
be checked by many efficient methods, e.g., eigenvalue method,
root locus method, value set method, etc.
\cite{ABK93,Bar94,BCK95}.

\section{Conclusions}

We have constructively proved that, for any two $n$-$th$ order
polynomials $a(s)$ and $b(s),$ the Hurwitz stability of their
convex combination is necessary and sufficient for the existence
of a polynomial $c(s)$ such that $c(s)/a(s)$ and $c(s)/b(s)$ are
both strictly positive real. By using similar method, we can also
constructively prove the existence of SPR synthesis for low order
($n\leq 4$) interval polynomials. Namely, when $n\leq 4,$ the
Hurwitz stability of the four Kharitonov vertex polynomials is
necessary and sufficient for the existence of a fixed polynomial
such that the ratio of this polynomial to any polynomial in the
interval polynomial set is SPR invariant
\cite{WY99,WY00,WY01a,Yu98,YW00,YW01}. The SPR synthesis problem
for high order interval polynomials is currently under
investigation.

 }}

\vskip 20pt
\vspace*{1\baselineskip}


\begin{thebibliography}{}

\bibitem{ABK93} {Ackermann, J., Bartlett, A., Kaesbauer, D., Sienel, W., and
Steinhauser, R., \textit{Robust Control: Systems with Uncertain
Physical Parameters.} Springer-Verlag, Berlin, 1993.}

\bibitem{ADK90} {Anderson, B. D. O., Dasgupta, S., Khargonekar, P., Kraus, F. J., and Mansou, M.,
Robust strict positive realness: characterization and
construction, \textit{IEEE Trans. Circuits Syst.}, 1990, CAS-37:
869-876.}

\bibitem{AM70} {Anderson, B. D. O. and Moore, J. B.,
\textit{Linear Optimal Control.} New York: Prentice Hall, 1970.}

\bibitem{Bar94} {Barmish, B. R.,
\textit{New Tools for Robustness of Linear Systems.} New York:
MacMillan Publishing Company, 1994.}

\bibitem{BZ93}  {Betser, A. and Zeheb, E., Design of robust strictly positive
real transfer functions. \textit{IEEE Trans. Circuits Syst.; Part
I,} 1993, \textbf{CAS-40:} 573-580.}

\bibitem{BCK95}  {Bhattacharyya, S. P., Chapellat, H., and Keel, L . H ., \textit{%
Robust Control: The Parametric Approach.} New York: Prentice Hall,
1995.}

\bibitem{Bia02}  {Bianchini, G., Synthesis of robust
strictly positive real discrete-time systems with $l_2$ parametric pertuebations. \textit{%
IEEE Trans. Circuits Syst.; Part I,} 2002, \textbf{CAS-49:}
1221-1225.}

\bibitem{BTV01}  {Bianchini, G., Tesi, A., and Vicino, A., Synthesis of robust
strictly positive real systems with $l_2$ parametric uncertainty. \textit{%
IEEE Trans. Circuits Syst.; Part I,} 2001, \textbf{CAS-48:}
438-450.}

\bibitem{CDB91}  {Chapellat, H., Dahleh, M., and Bhattacharyya, S. P., On robust
nonlinear stability of interval control systems. \textit{IEEE
Trans. Automat. Contr.,} 1991, \textbf{AC-36:} 59-69.}

\bibitem{DB87}  {Dasgupta, S. and Bhagwat, A. S., Conditions for designing
strictly positive real transfer functions for adaptive output
error identification. \textit{IEEE Trans. Circuits Syst.,} 1987,
\textbf{CAS-34:} 731-737.}

\bibitem{DV75}  {Desoer, C. A., Vidyasagar, M. \textit{Feedback Systems:
Input-Output Properties.} San Diego: Academic Press, 1975.}

\bibitem{Hen02}  {Henrion, D., Linear matrix inequalities for robust
strictly positive real design. \textit{%
IEEE Trans. Circuits Syst.; Part I,} 2002, \textbf{CAS-49:}
1017-1020.}

\bibitem{HHX89}  {Hollot, C. V., Huang, L., and Xu, Z. L., Designing strictly
positive real transfer function families: A necessary and
sufficient condition for low degree and structured families.
\textit{Proc. of Mathematical Theory of Network and Systems,}
(eds. Kaashoek, M. A., Van Schuppen, J. H., Ran, A. C. M.),
Boston, Basel, Berlin: Birkh\"auser, 1989, 215-227.}

\bibitem{HHX90}  {Huang, L., Hollot, C. V., Xu, Z. L., Robust analysis of
strictly positive real function set. Preprints of \textit{The
Second Japan-China Joint Symposium on Systems Control Theory and
its Applications,} 1990, 210-220.}

\bibitem{Kal63}  {Kalman, R. E., Lyapunov functions for the problem of Lur'e in
automatic control. \textit{Proc. Nat. Acad. Sci.}(USA), 1963,
\textbf{49:} 201-205.}

\bibitem{Lan79}  {Landau, I. D., \textit{Adaptive Control: The Model Reference
Approach.} New York: Marcel Dekker, 1979.}

\bibitem{MA98}  {Marquez, H. J. and Agathoklis, P., On the existence of robust
strictly positive real rational functions. \textit{IEEE Trans.
Circuits Syst.; Part I,} 1998, \textbf{CAS-45:} 962-967.}

\bibitem{MP01}  {Mosquera, C. and Perez, F., Algebraic solution to the robust
SPR problem for two polynomials. \textit{Automatica,} 2001,
\textbf{37:} 757-762.}

\bibitem{PD97}  {Patel, V. V. and Datta, K. B., Classification of units in $%
H_\infty $ and an alternative proof of Kharitonov's theorem.
\textit{IEEE Trans. Circuits Syst.; Part I,} 1997,
\textbf{CAS-44:} 454-458.}

\bibitem{Pop73}  {Popov, V. M., \textit{Hyperstability of Control Systems.} New
York: Springer-Verlag, 1973.}

\bibitem{WH91}  {Wang, L. and Huang, L., Finite verification of strict
positive realness of interval rational functions. \textit{Chinese
Science Bulletin, }1991, \textbf{36:} 262-264.}

\bibitem{WY99}  {Wang, L. and Yu, W. S., A new approach to robust synthesis of
strictly positive real transfer functions. \textit{Stability and
Control: Theory and Applications,} 1999, \textbf{2:} 13-24.}

\bibitem{WY00}  {Wang, L. and Yu, W. S., Complete characterization of strictly
positive real regions and robust strictly positive real synthesis
method. \textit{Science in China,} 2000, \textbf{(E)43:}\ 97-112.}

\bibitem{WY01}  {Wang, L. and Yu, W. S., On robust stability of polynomials
and robust strict positive realness of transfer functions.
\textit{IEEE Trans. on Circuits and Syst.; Part I,} 2001,
\textbf{CAS-48:} 127-128.}

\bibitem{WY01a}  {Wang, L. and Yu, W. S., Robust SPR synthesis for low-order
polynomial segments and interval polynomials. \textit{Proceedings
of the American Control Conference (ACC 2001),} Crystal Gateway
Marriott, Arlington, VA, USA, 2001, 3612-3617.}

\bibitem{XWY02}  {Xie, L. J., Wang, L., and Yu, W. S., A new geometric
algorithm with order reduction for robust strictly positive real
synthesis. \textit{2002 The 41st IEEE Conference on Decision and
Control (CDC 2002),} Las Vegas, NV, USA, 2002.}

\bibitem{Yu98}  {Yu, W. S., \textit{Robust Strictly Positive Real Synthesis
and Robust Stability Analysis.} PhD Thesis, Peking University,
Beijing, 1998.}

\bibitem{YH99}  {Yu, W. S. and Huang, L., A necessary and sufficient conditions
on robust SPR stabilization for low degree systems.
\textit{Chinese Science Bulletin,} 1999, \textbf{44:} 517-520.}

\bibitem{YW99}  {Yu, W. S. and Wang, L., Some remarks on the definition of
strict positive realness of transfer Functions.
\textit{Proceedings of Chinese Conference on Decision and
Control,} Northeast University Press, Shenyang, 1999, 135-139.}

\bibitem{YW00}  {Yu, W. S. and Wang, L., Design of strictly positive real
transfer functions. \textit{IFAC Symposium on Computer Aided
Control Systems Desgin (CACSD 2000),} Salford, UK, 2000.}

\bibitem{YW01}  {Yu, W. S. and Wang, L., Anderson's claim on fourth-order SPR
synthesis is true. \textit{IEEE Trans. Circuits Syst.; Part I,}
2001, \textbf{CAS-48:} 506-509.}

\bibitem{YW01a}  {Yu, W. S. and Wang, L., Robust SPR synthesis for fourth-order
convex combinations. \textit{Progress in natural science,} 2001,
\textbf{11:} 461-467.}

\bibitem{YW01b}  {Yu, W. S. and Wang, L., Robust strictly positive real
synthesis for convex combination of the fifth-order polynomials. \textit{%
Proceedings of the IEEE Symposium on Circuits and Systems
Conference (ISCAS 2001),} Sydney, Austrlia, 2001, Volume 1:
739-742.}

\bibitem{YW03}  {Yu, W. S. and Wang, L., Robust strictly positive real
synthesis for convex combination of the sixth-order polynomials. \textit{%
Proceedings of the American Control Conference (ACC 2003),} USA,
2003 (submitted).}

\bibitem{YWT99}  {Yu, W. S., Wang, L., and Tan, M., Complete characterization of
strictly positive realness regions in coefficient space.
\textit{Proceedings of the IEEE Hong Kong Symposium on Robotics
and Control,} Hong Kong, 1999, 259-264.}

\end{thebibliography}
\end{document}